\documentclass[12pt,a4paper]{amsart}
\usepackage{t1enc}
\usepackage[latin1]{inputenc}
\usepackage[english]{babel}
\usepackage{hyperref}
\usepackage{graphicx}

\usepackage{latexsym}
\usepackage{amsmath,amsfonts,amssymb,amsthm}
\usepackage{mathptmx}

\begin{document}

\title {
	Adapted complex tubes on the symplectization of pseudo-Hermitian manifolds
}

\author{Giuseppe Tomassini}
\address{
	G. Tomassini:
	Scuola Normale Superiore,
	Piazza dei Cavalieri,
	7 --- I-56126 Pisa, Italy}
\email{g.tomassini@sns.it}

\author{Sergio Venturini}
\address{
	S. Venturini:
	Dipartimento Di Matematica,
	Universit\`{a} di Bologna,
	\,\,Piazza di Porta S. Donato 5 ---I-40127 Bologna,
	Italy}
\email{venturin@dm.unibo.it}


%/intest
%\magnification=\magstep1
%\font\chapt=cmbx10 scaled 1700
%\font\sc=cmcsc10
%\font\ninerm=cmr8
%\font\nineit=cmmi8
%\font\ninesl=cmsl8
%\font\titolo=cmbx9 scaled \magstep 3
%
%
%numeri....
\def\R{{\rm I\kern-.185em R}}
\def\C{{\rm\kern.37em\vrule height1.4ex width.05em depth-.011em\kern-.37em C}}
\def\N{{\rm I\kern-.185em N}}
\def\Z{{\bf Z}}
\def\Q{{\mathchoice{\hbox{\rm\kern.37em\vrule height1.4ex width.05em 
depth-.011em\kern-.37em Q}}{\hbox{\rm\kern.37em\vrule height1.4ex width.05em 
depth-.011em\kern-.37em Q}}{\hbox{\sevenrm\kern.37em\vrule height1.3ex 
width.05em depth-.02em\kern-.3em Q}}{\hbox{\sevenrm\kern.37em\vrule height1.3ex
width.05em depth-.02em\kern-.3em Q}}}}
\def\P{{\rm I\kern-.185em P}}
\def\H{{\rm I\kern-.185em H}}
%
%insiemi
\def\Aleph{\aleph_0}
\def\ALEPH#1{\aleph_{#1}}
\def\sset{\subset}\def\ssset{\sset\sset}
%
%funzioni
\def\bar#1{\overline{#1}}
\def\dim{\mathop{\rm dim}\nolimits}
\def\half{\textstyle{1\over2}}
\def\Half{\displaystyle{1\over2}}
\def\mlog{\mathop{\half\log}\nolimits}
\def\Mlog{\mathop{\Half\log}\nolimits}
\def\Det{\mathop{\rm Det}\nolimits}
\def\Hol{\mathop{\rm Hol}\nolimits}
\def\Aut{\mathop{\rm Aut}\nolimits}
\def\Re{\mathop{\rm Re}\nolimits}
\def\Im{\mathop{\rm Im}\nolimits}
\def\Ker{\mathop{\rm Ker}\nolimits}
\def\Fix{\mathop{\rm Fix}\nolimits}
\def\Res{\mathop{\rm Res}\nolimits}
\def\sp{\mathop{\rm sp}\nolimits}
\def\id{\mathop{\rm id}\nolimits}
\def\Trace{\mathop{\rm Tr}\nolimits}
\def\cancel#1#2{\ooalign{$\hfil#1/\hfil$\crcr$#1#2$}}
\def\prevoid{\mathrel{\scriptstyle\bigcirc}}
\def\void{\mathord{\mathpalette\cancel{\mathrel{\scriptstyle\bigcirc}}}}
\def\n{{}|{}\!{}|{}\!{}|{}}
\def\abs#1{\left|#1\right|}
\def\norm#1{\left|\!\left|#1\right|\!\right|}
\def\nnorm#1{\left|\!\left|\!\left|#1\right|\!\right|\!\right|}
%
%integrali superiore ed inferiore
\def\upperint{\int^{{\displaystyle{}^*}}}
\def\lowerint{\int_{{\displaystyle{}_*}}}
\def\Upperint#1#2{\int_{#1}^{{\displaystyle{}^*}#2}}
\def\Lowerint#1#2{\int_{{\displaystyle{}_*}#1}^{#2}}
%
%altro
\def\rem #1::#2\par{\medbreak\noindent{\bf #1}\ #2\medbreak}
\def\proclaim #1::#2\par{\removelastskip\medskip\goodbreak{\bf#1:}
\ {\sl#2}\medskip\goodbreak}
\def\ass#1{{\rm(\rmnum#1)}}
\def\assertion #1:{\Acapo\llap{$(\rmnum#1)$}$\,$}
\def\Assertion #1:{\Acapo\llap{(#1)$\,$}}
\def\acapo{\hfill\break\noindent}
\def\Acapo{\hfill\break\indent}
\def\proof{\removelastskip\par\medskip\goodbreak\noindent{\it Proof.\/\ }}
\def\prova{\removelastskip\par\medskip\goodbreak
\noindent{\it Dimostrazione.\/\ }}
\def\risoluzione{\removelastskip\par\medskip\goodbreak
\noindent{\it Risoluzione.\/\ }}
\def\qed{{$\Box$}\par\smallskip}
\def\BeginItalic#1{\removelastskip\par\medskip\goodbreak
\noindent{\it #1.\/\ }}
\def\iff{if, and only if,\ }
\def\sse{se, e solo se,\ }
\def\rmnum#1{\romannumeral#1{}}
\def\Rmnum#1{\uppercase\expandafter{\romannumeral#1}{}}
\def\smallfrac#1/#2{\leavevmode\kern.1em
\raise.5ex\hbox{\the\scriptfont0 #1}\kern-.1em
/\kern-.15em\lower.25ex\hbox{\the\scriptfont0 #2}}
%
%delimitatori
\def\Left#1{\left#1\left.}
\def\Right#1{\right.^{\llap{\sevenrm
\phantom{*}}}_{\llap{\sevenrm\phantom{*}}}\right#1}
%
%
%bibliografia
%\def\dimens{3em}
%\newcount\qqrefno
%\qqrefno=0
%\def\qqrefnum{\global\advance\qqrefno by 1
%\noindent\rlap{[\number\qqrefno]}\hbox to \dimens{}\hangindent=\dimens}
%\def\references{\bigskip\noindent{\bf References.}\bigskip}
%\def\art #1 : #2 ; #3 ; #4 ; #5 ; #6. \par{\qqrefnum
%#1, {\sl#2}, #3, {\bf#4}, (#5), #6.\par\smallskip}
%\def\book #1 : #2 ; #3 ; #4. \par{\qqrefnum#1, {\bf#2}, #3, #4.\par\smallskip}
%\def\freeart #1 : #2 ; #3. \par{\qqrefnum#1, {\sl#2}, #3.\par\smallskip}
%
%bibliografiabis
\def\dimens{3em}
\def\symb[#1]{\noindent\rlap{[#1]}\hbox to \dimens{}\hangindent=\dimens}
\def\references{\bigskip\noindent{\bf References.}\bigskip}
\def\art #1 : #2 ; #3 ; #4 ; #5 ; #6. \par{#1, 
{\sl#2}, #3, {\bf#4}, (#5), #6.\par\smallskip}
\def\book #1 : #2 ; #3 ; #4. \par{#1, {\bf#2}, #3, #4.\par\smallskip}
\def\freeart #1 : #2 ; #3. \par{#1, {\sl#2}, #3.\par\smallskip}
%
%
%bibliografia per riviste tedesche
%\def\dimens{3em}
%\newcount\qqrefno
%\qqrefno=0
%\def\qqrefnum{\global\advance\qqrefno by 1
%\noindent\hbox to \dimens{}\hangindent=\dimens
%\llap{\number\qqrefno.\ \ \ }}
%\def\references{\bigskip\noindent{\bf References.}\bigskip}
%\def\art #1 : #2 ; #3 ; #4 ; #5 ; #6. \par{\qqrefnum
%#1: #2. #3 {\bf#4}, #6, (#5)\par\smallskip}
%\def\book #1 : #2 ; #3. \par{\qqrefnum#1: #2. #3.\par\smallskip}
%\def\freeart #1 : #2 ; #3. \par{\qqrefnum#1: #2. #3.\par\smallskip}
%
%
%nome e indirizzi.
\def\name{\hbox{Sergio Venturini}}
\def\snsaddress{\indent
\vbox{\bigskip\bigskip\bigskip
\name
\hbox{Scuola Normale Superiore}
\hbox{Piazza dei Cavalieri, 7}
\hbox{56126 Pisa (ITALY)}
\hbox{FAX 050/563513}}}
\def\cassinoaddress{\indent
\vbox{\bigskip\bigskip\bigskip
\name
\hbox{Universit\`a di Cassino}
\hbox{via Zamosch 43}
\hbox{03043 Cassino (FR)}
\hbox{ITALY}}}
\def\bolognaaddress{\indent
\vbox{\bigskip\bigskip\bigskip
\name
\hbox{Dipartimento di Matematica}
\hbox{Universit\`a di Bologna}
\hbox{Piazza di Porta S. Donato 5}
\hbox{40127 Bologna (BO)}
\hbox{ITALY}
\hbox{venturin@dm.unibo.it}
}}
\def\homeaddress{\indent
\vbox{\bigskip\bigskip\bigskip
\name
\hbox{via Garibaldi, 7}
\hbox{56124 Pisa (ITALY)}}}
\def\doubleaddress{
\vbox{
\hbox{\name}
\hbox{Universit\`a di Cassino}
\hbox{via Zamosch 43}
\hbox{03043 Cassino (FR)}
\hbox{ITALY}
\smallskip
\hbox{and}
\smallskip
\hbox{Scuola Normale Superiore}
\hbox{Piazza dei Cavalieri, 7}
\hbox{56126 Pisa (ITALY)}
\hbox{FAX 050/563513}}}
\def\sergio{{\rm\bigskip
\centerline{Sergio Venturini}
\leftline{\bolognaaddress}
\bigskip}}
%
%
%allineamenti formule matematiche
%
%\centeredeq{...}
%                      |*  [ ]&*   [ ]&....|
%                      ..................
%\leftcenteredeq
%|*  [ ]&*   [ ]&....|
%..................
%
%\def\centeredeq#1{\vcenter{\halign
%{&$\strut\displaystyle{##}\hfil\strut\qquad$\cr#1}}
%\def\leftcenteredeq#1{\halign
%{&$\strut\displaystyle{##}\hfil\strut\qquad$\cr#1}}
%
%@@@@@@@
%%\input{../greek}
%macrosl

\newtheorem{theorem}{Theorem}[section]
\newtheorem{proposition}{Proposition}[section]
\newtheorem{lemma}{Lemma}[section]
\newtheorem{corollary}{Corollary}[section]
\newtheorem{remark}{Remark}[section]
\newtheorem{definition}{Definition}[section]

\newtheorem{teorema}{Teorema}[section]
\newtheorem{proposizione}{Proposizione}[section]
\newtheorem{corollario}{Corollario}[section]
\newtheorem{osservazione}{Osservazione}[section]
\newtheorem{definizione}{Definizione}[section]
\newtheorem{esempio}{Esempio}[section]
\newtheorem{esercizio}{Esercizio}[section]
\newtheorem{congettura}{Congettura}[section]

\bibliographystyle{abbrv}

\def\RManBase{M}
\def\RFormBase{\theta}
\def\EqBase{E}
\def\ReebVField{{\xi_0^\RFormBase}}
\def\ReebXVField{{\xi^\RFormBase}}
\def\domJ{U}
\def\fvar{\alpha}
\def\holpar{\phi}
\def\DiffeomA{\Phi}
\def\DiffeomB{\Psi}
\def\CMap{\Psi^c}
\def\CManBase{{\tilde{\RManBase}}}
\def\CEmbedding{j}
\def\NeighbrA{V}
\def\NeighbrB{W}
\def\NeighbrAA{{V_1}}
\def\NeighbrBB{{V_2}}
\def\CR{{\rm CR}}
\def\dd{{\rm d}}
\def\dc{{\dd^c}}
\def\res{\mathop{\hbox{\vrule height 7pt width .5pt depth 0pt\vrule height .5pt width 6pt depth 0pt}}\nolimits}

\begin{abstract}
Let $(\RManBase,\RFormBase)$ be a pseudo-Hermitian space
of real dimension $2n+1$, that is $\RManBase$ is a $\CR-$manifold of dimension
$2n+1$ and $\RFormBase$ is a contact form on $\RManBase$
giving the Levi distribution $HT(\RManBase)\sset T\RManBase$.
Let $\RManBase^\RFormBase\sset T^*\RManBase$ be the canonical
symplectization of $(\RManBase,\RFormBase)$ and $\RManBase$
be identified with the zero section of $\RManBase^\RFormBase$.
Then $\RManBase^\RFormBase$ is a manifold of real dimension $2(n+1)$
which admit a canonical foliation by surfaces parametrized
by $\mathbb{C}\ni t+i\sigma\mapsto \holpar_p(t+i\sigma)=\sigma\RFormBase_{g_t(p)}$,
where $p\in\RManBase$ is arbitrary and $g_t$
is the flow generated by the Reeb vector field associated to
the contact form $\RFormBase$.

Let $J$ be an (integrable) complex structure defined in a neighbourhood
$\domJ$ of $\RManBase$ in $\RManBase^\RFormBase$.
We say that the pair $(\domJ,J)$ is an {adapted complex tube}
on $\RManBase^\RFormBase$ if all the parametrizations
$\holpar_p(t+i\sigma)$ defined above are holomorphic on $\holpar_p^{-1}(\domJ)$.

In this paper we prove that if $(\domJ,J)$ is an adapted complex tube on $\RManBase^\RFormBase$,
then the real function $\EqBase$ on $\RManBase^\RFormBase\sset T^*\RManBase$
defined by the condition
$\fvar=\EqBase(\fvar)\RFormBase_{\pi(\fvar)}$,
for each $\fvar\in\RManBase^\RFormBase$,
is a canonical equation for $\RManBase$
which satisfies the homogeneous Monge-Amp\`ere equation
$(\dd\dc\EqBase)^{n+1}=0$.

We also prove that if $\RManBase$ and $\RFormBase$ are real analytic
then the symplectization $\RManBase^\RFormBase$ admits an unique
maximal adapted complex tube.

\end{abstract}

\maketitle
%\tableofcontents

%\nocite{book:DragominTomassini}
%\nocite{book:KobayashiNomizuA}
%\nocite{book:KobayashiNomizuB}
%\nocite{book:KushnerLychaginRubtsov}
%\nocite{book:GilbertBuchanan}
%\nocite{book:BlairContactEtc}
%\nocite{article:AndreottiFredricksEmbeddingCR}
%\nocite{article:BedfordKalka}
%\nocite{article:DuchampKalka}
%\nocite{article:BedfordBurns}
%\nocite{article:TV2ArxivContactMA}

%Proposition 1.5 pag. 387 of \cite{article:BedfordBurns}: problema di Cauchy.

\section{\label{section:A}Adapted complex tubes}

In this paper we apply some results obtained by the authors in
\cite{article:TV2ArxivContactMA} to the study of
suitable ``adapted'' integrable complex structure on open neighbourhood 
of a pseudo-Hermitian manifold in its symplectization.

We follow \cite{book:DragominTomassini} for standard notations
in complex and $\CR-$geometry.

Given a smooth real differentiable manifold $\RManBase$ of dimension $m$
we denote by $T^*\RManBase$ its cotangent bundle and
$\pi:T^*\RManBase\to\RManBase$ the canonical projection.

Let $\RFormBase$ be a smooth $1-$form on $M$
such that $\RFormBase_p\neq0$ for each $p\in\RManBase$.
Then the subset of by $T^*\RManBase$
$$
	\RManBase^\RFormBase=\bigl\{\fvar\in T^*\RManBase\mid
		\fvar\wedge\RFormBase_{\pi(\fvar)}=0\bigr\},
$$
is the line bundle on $\RManBase\sset T^*\RManBase$ of the forms
$\fvar$ which are proportional to the form $\RFormBase_{\pi(\fvar)}$.
We denote by $\EqBase:\RManBase^\RFormBase\to\mathbb{R}$ the function
which satisfies
$$
	\fvar=\EqBase(\fvar)\RFormBase_{\pi(\fvar)}
$$
for each $\fvar\in T^*\RManBase$.

Denoting by $\RManBase_0$ the zero section of $T^*\RManBase$
we see that for each $\fvar\in\RManBase^\RFormBase$ we have
$\EqBase(\fvar)=0$ if, and only if, $\fvar\in\RManBase_0$.

We freely identify the manifold $\RManBase$ with the zero section
$\RManBase_0$ in $\RManBase^\RFormBase$.

It easy to show that $\RManBase^\RFormBase$ is a smooth closed submanifold of
$T^*\RManBase$ of (real) dimension $m+1$ and
$\EqBase$ is a smooth function on $\RManBase$.

We also denote by $X^\RFormBase$
the vector field on $T^*\RManBase$ which is the infinitesimal
generator of the one parameter group of transformation defined by
$$h_t^\RFormBase(\fvar)=\fvar+t\RFormBase_{\pi(\fvar)}$$
for each $\fvar\in T^*\RManBase$ and $t\in\mathbb{R}$.

Let $x=(x_1,\ldots,x_m)$ be a local coordinate system on $\RManBase$ and let
$(x,p)=(x_1,\ldots,x_m,p_1,\ldots,p_m)$ the corresponding local
coordinate system on $T^*\RManBase$, so that
$$\RFormBase=\sum_{i=1}^mp_i(\RFormBase)dx_i.$$
Then we have
$$X^\RFormBase=\sum_{i=1}^mp_i(\RFormBase)\frac{\partial}{\partial p_i}.$$
Clearly the vector field $X^\RFormBase$ is tangent to $\RManBase^\RFormBase$
and $X^\RFormBase(\EqBase)=1$.
Thus the smooth function $\EqBase$ has no critical points on $\RManBase^\RFormBase$.

Let now $(\RManBase,\RFormBase)$ be a \emph{pseudo-Hermitian manifold},
that is $\RManBase$ is a orientable $\CR-$manifold of dimension $2n+1$
with non degenerate Levi form and $\RFormBase$ is a non degenerate real $1-$form which vanishes on the Levi distribution $HT(\RManBase)\sset T\RManBase$ and defines the \emph{pseudo-Hermitian structure} on $\RManBase$. Then $(\RManBase,\RFormBase)$ is a contact manifold with volume form $\RFormBase\wedge(\dd\RFormBase)^{n}$ and there exists a vector field $\ReebVField$, the \emph{Reeb vector field}, which is characterized by the conditions
$\ReebVField\res\RFormBase=1$ and $\ReebVField\res\dd\RFormBase=0$. The pull-back of $\dd\RFormBase$ is a symplectic form on $\RManBase^\RFormBase$ that will be called the \emph{symplectization} of $\RManBase$.

Let $g_t$ be the one parameter (local) group of transformations
associated to $\ReebVField$, that is for each $p\in\RManBase$ the map
$$
g_{\cdot}(p):I_p\to\RManBase,\quad t\mapsto g_t(p)\RManBase,\quad I_p\sset\mathbb{R}
$$
is the maximal integral curve of $\ReebVField$ such that $g_0(p)=p$.

For each $p\in\RManBase$ let $\tilde{I}_p=\{z=t+i\sigma\in\mathbb{C}\mid t\in I_p\}$ and the map
$$
\holpar_p=\holpar_p^\RFormBase:\tilde{I}_p\to\RManBase^\RFormBase
$$
defined by the formula
$$
\holpar_p^\RFormBase(z)=\holpar^\RFormBase(t+i\sigma)=\sigma\RFormBase_{g_t(p)}.
$$

We also denote by $\ReebXVField$ the vector field on $\RManBase^\RFormBase$
which is the infinitesimal generator of the one parameter local group of transformations
defined by
$$
g_t^\RFormBase(\fvar)=\EqBase(\fvar)\RFormBase_{g_t(p)},
$$
where $p=\pi(\fvar)$.

Then the vector field $\ReebXVField$ extends $\ReebVField$ and
it is easy to show that the one parameter groups $g_t^\RFormBase$
and $h_t^\RFormBase$ commute and generate a two dimensional distribution on
$\RManBase^\RFormBase$ whose leaves are parametrized by the maps
$\holpar_p^\RFormBase$.
Moreover we have
\begin{equation}\label{eq::CFA}
\left.
		\begin{array}{l}
			\ReebXVField(\EqBase)=0,\\
			X^\RFormBase(\EqBase)=1,\\
			\left[\ReebXVField,X^\RFormBase\right]=0.
		\end{array}
\right.
\end{equation}

We now give the main definition of this paper.

\begin{definition}\label{def::AdaptedContactTube}
Let $(\RManBase,\RFormBase)$ be a \emph{pseudo-Hermitian manifold}
of real dimension $2n+1$.
Let the manifold $\RManBase$ be identified with the zero section
$\RManBase_0$ in $\RManBase^\RFormBase$.

An {adapted complex tube} on $\RManBase^\RFormBase$ is a pair
$(\domJ,J)$ where $\domJ$ is a open subset of $\RManBase^\RFormBase$
containing the zero section $\RManBase_0$ and $J$ is an (integrable)
complex structure on $\domJ$ which satisfies the following conditions:
\begin{enumerate}
	\item for each $p\in\RManBase$ the set $U_p=\{\fvar\in U\mid \pi(\fvar)=p\}$
	is connected;
	\item the restriction of $J$ to the Levi distribution $HT(\RManBase)$ is the given
	$\CR$ structure of $HT(\RManBase)$;
	\item if $p\in\RManBase$ the restriction of the map $\holpar_p^\RFormBase$ to
	the open set $\tilde{I}_p\cap({\holpar_p^\RFormBase})^{-1}(\domJ)$ is
	holomorphic.
\end{enumerate}
An adapted complex tube $(\domJ,J)$ with $J$ of class $C^k$, $k=1$,$\ldots$,$\infty,\omega$, is said to be \emph{$C^k-$maximal},
 if for any other adapted complex tube on $\RManBase^\RFormBase$
such that $J'$ is of class $C^k$, $\domJ\sset\domJ'$ and $J'_{\vert U}=J$ one has $\domJ=\domJ'$ (and hence $J'=J$).
\end{definition}
Adapted tubes of class $C^\infty$ or $C^\omega$ are said to be \emph{smooth} or  \emph{real analytic}, respectively.

The main results of this paper are the following.

\begin{theorem}\label{thm::AdaptedMA}
Let $(\RManBase,\RFormBase)$ be a pseudo-Hermitian manifold
of real dimension $2n+1$ with  $\RManBase$ and $\RFormBase$ of class $C^\infty$.
Let $(\domJ,J)$ be any smooth adapted complex tube on $\RManBase^\RFormBase$.
Then the function $E$ satisfies
	\begin{equation}\label{eq::CauchyMongeAmpere}
		\left\{
		\begin{array}{l}
		(\dd\dc\EqBase)^{n+1}=0\ {\rm on}\ \domJ,\\
		\\
		\dd\EqBase\wedge\dc\EqBase\wedge(\dd\dc\EqBase)^{n}\neq0\ {\rm near}\ \RManBase,\\
		\\
		\EqBase_{|\RManBase}=0,\\
		\\
		\dc\EqBase_{|T(\RManBase)}=-\RFormBase.
		\end{array}
		\right.
	\end{equation}
\end{theorem}

\begin{theorem}\label{thm::AdaptedRealAnalytic}
Let $(\RManBase,\RFormBase)$ be a pseudo-Hermitian manifold
of real dimension $2n+1$ with  $\RManBase$ and $\RFormBase$ real analytic.
Then there exists a unique real analytic maximal adapted complex tube
$(\domJ,J)$ on $\RManBase^\RFormBase$.
\end{theorem}

\medskip
{\bf Remark 1.1. }
Adapted complex structure on the tangent (and cotangent) bundle
a Riemannian manifols are studied in
\cite{article:LempertSzoke},
\cite{article:Szoke} and, independentely, in
\cite{article:GuilleminStenzelI}.

\section{\label{section:Proof:A}Proof of Theorem 1.1}
Let $(\RManBase,\RFormBase)$ be a fixed pseudo-Hermitian manifold.

\begin{lemma}\label{lemma::AdaptedJ}
If $(\domJ,J)$ is any smooth adapted complex tube
on $\RManBase^\RFormBase$ then

\begin{equation}\label{eq::AdaptedJ}
	J(\ReebXVField)=X^\RFormBase.
\end{equation}
\end{lemma}

\proof
Let $\fvar\in\domJ\sset\RManBase^\RFormBase$ and $p=\pi(\fvar)$.
Then  $\fvar=\sigma_0\RFormBase_p$, where $\sigma_0=\EqBase(\fvar)$.
By definition of $\ReebXVField$
\begin{eqnarray*}
	\ReebXVField(\fvar)&=&\left.\frac{d}{dt}g_t^\RFormBase(\fvar)\right|_{t=0}
		=\left.\frac{d}{dt}\bigl(\sigma_0\RFormBase_{g_t(p)}\bigr)\right|_{t=0}\cr
		&=&\left.\frac{d}{dt}\holpar_p(t+i\sigma_0)\right|_{t=0}.
\end{eqnarray*}
Since the map $\holpar_p(t+i\sigma)$ is holomorphic we obtain
\begin{eqnarray*}
	J\bigl(\ReebXVField(\fvar)\bigr)
		&=&\left.\frac{d}{d\sigma}\holpar_p\bigl(i(\sigma_0+\sigma)\bigr)\right|_{\sigma=0}
		=\left.\frac{d}{d\sigma}\bigl(\fvar+\sigma\RFormBase_p\bigr)\right|_{\sigma=0}\cr
		&=&\left.\frac{d}{d\sigma}h_\sigma^\RFormBase(\fvar)\right|_{\sigma=0}
		=X^\RFormBase(\fvar),
\end{eqnarray*}
as desired.

\qed

Using the terminology of \cite{article:TV2ArxivContactMA}
we observe that by (\ref{eq::CFA}) and (\ref{eq::AdaptedJ})
the pair $(\ReebXVField,-\EqBase)$ is a smooth
\emph{calibrated foliation}, that is

\begin{eqnarray*}
	&&[\ReebXVField,J(\ReebXVField)]=0,\cr
	&&\dd\EqBase(\ReebXVField)=0,\cr
	&&\dc\EqBase(\ReebXVField)=-1.\cr
\end{eqnarray*}
Moreover the function $\EqBase$ vanishes on $\RManBase$ and 
the vector field $\ReebXVField$ extends
the Reeb vector field $\ReebVField$ on $\RManBase$. Since the Reeb vector field $\ReebVField$ is an infinitesimal
symmetry for the contact distribution $HT(\RManBase)$
from  of \cite[Theorem 3.1]{article:TV2ArxivContactMA} it follows 
that the form $\ReebXVField\res\dd\dc\EqBase=L_\ReebXVField\dc\EqBase$
vanishes identically on $\domJ$.
Here $L_\ReebXVField\dc\EqBase$ stands for the Lie
derivative of the form $\dc\EqBase$
with respect to the vector field $\ReebXVField$.
Since $\ReebXVField\neq0$ on $\domJ$ it follows that
the form $\dd\dc\EqBase$ satisfies the Monge-Amp\`ere equation
$(\dd\dc\EqBase)^{n+1}=0$ on $\domJ$.

Clearly, the function $\EqBase$ vanishes exactly on $\RManBase$
and hence the restriction of $\dc\EqBase$ to $T\RManBase$ is
$\lambda\RFormBase$ for some (smooth) function $\lambda:\RManBase\to\mathbb{R}$.
For each $p\in\RManBase$ we have
\begin{eqnarray*}
	\lambda(p)&=&\lambda(p)\RFormBase_p\bigl(\ReebVField(p)\bigr)
		=\dc\EqBase\bigl(\ReebVField(p)\bigr)\cr
	&=&\dc\EqBase\bigl(\ReebXVField(p)\bigr)=
		-\dd\EqBase\bigl(X^\RFormBase(p))\bigr)=-1
\end{eqnarray*}
and hence $\dc\EqBase_{|T(\RManBase)}=-\RFormBase$.

Since $\RFormBase\wedge(\dd\RFormBase)^n$ is a volume form on $\RManBase$ we have
$\dd\EqBase\wedge\dc\EqBase\wedge(\dd\dc\EqBase)^{n}\neq0$
in a neighbourhood of $\RManBase$.

This completes the proof of Theorem \ref{thm::AdaptedMA} .

\section{\label{section:Proof:B}Proof of Theorem 1.2}%\ref{thm::AdaptedRealAnalytic}
Let $(\RManBase,\RFormBase)$ be a fixed pseudo-Hermitian manifold.

Let $\DiffeomA:\RManBase\times\mathbb{R}\to\RManBase^\RFormBase$ be the map
defined by $\DiffeomA(p,\sigma)=\sigma\RFormBase_p$ for each $p\in\RManBase$
and each $\sigma\in\mathbb{R}$.
Then $\DiffeomA$ is a diffeomorphism between $\RManBase\times\mathbb{R}$
and $\RManBase^\RFormBase$.

By \cite{article:AndreottiFredricksEmbeddingCR} there exists a pair $(\CManBase, \CEmbedding)$
where $\CManBase$ is a complex manifold and $\CEmbedding$ is a real analytic embedding of
$\RManBase$ in $\CManBase$.
Then, by analytic extension, there exist an open neighbourhood $V$ of $\RManBase\times\{0\}$ in
$\RManBase\times\mathbb{C}$ and a map $\CMap:\NeighbrA\to\CManBase$ satisfying
\begin{enumerate}
	\item for each $p\in\RManBase$ the set $I_p=\{t\in\mathbb{R}\mid(p,t)\in\NeighbrA$
	is an open interval in $\mathbb{R}$ and $\CMap(p,t)=g_t(p)$ is the maximal integral curve
	of the vector field $\ReebVField$ such that $g_0(p)=0$, i.e. $\CMap(p,0)=0$ for each
	$p\in\RManBase$;
	\item for each $(p,t)\in\NeighbrA\cap\RManBase\times\mathbb{R}$ the set
	$\NeighbrA_{(p,t)}=\{\sigma\in\mathbb{R}\mid(p,t+i\sigma)\in\NeighbrA\}$
	is an open interval in $\mathbb{R}$ containing the origin $0\in\mathbb{R}$
	\item for each $p\in\RManBase$ the map $z=(t+i\sigma)\mapsto\CMap(p,z)$ is
	holomorphic when defined.
\end{enumerate}

Then, set $\NeighbrB=\{(p,\sigma)\in\RManBase\times\mathbb{R}\mid(p,i\sigma)\in\NeighbrA\}$
and define $\DiffeomB:\NeighbrB\to\CManBase$ by the formula
$\DiffeomB(p,\sigma)=\CMap(p,i\sigma)$.
Since $X^\RFormBase(p)=J\bigl(\ReebVField(p)\bigr)$ is transversal to $\RManBase$ for
each $p\in\RManBase$ it follows that,
shrinking $\NeighbrA$ (and hence $\NeighbrB$) and $\CManBase$ if necessary,
the map $\DiffeomB$ is a real analytic diffeomorphism between $\NeighbrB$ and $\CManBase$.

Let $\domJ=\DiffeomA(\NeighbrA)$ and $J$ be the pullback under
$\DiffeomA\circ\DiffeomB^{-1}$ of the complex structure of $\CManBase$. It is then easy to show that $(\domJ,J)$ is a
real analytic adapted complex tube on $\RManBase^\RFormBase$.

This proves the existence of real analytic adapted complex tubes.

Let now $(\domJ_1,J_1)$ and $(\domJ_2,J_2)$ be
two real analytic adapted complex tubes on $\RManBase^\RFormBase$.
We claim that $J_1$ and $J_2$ agree on the intersection $\domJ_1\cap\domJ_2$.
By an analytic continuation argument it suffices to prove that $J_1$ and $J_2$ agree
on a neighbourhood of $\RManBase$ in $\domJ_1\cap\domJ_2$.

Again by \cite{article:AndreottiFredricksEmbeddingCR}, there exist open neighbourhoods $\NeighbrAA\sset\domJ_1$ and $\NeighbrBB\sset\domJ_2$
of $\RManBase$ and a real analytic diffeomorphism $F:\NeighbrAA\to\NeighbrBB$ which is
the identity map on $\RManBase$ and which is holomorphic when $\NeighbrAA$ and $\NeighbrBB$
are endowed by the complex structures $J_1$ and $J_2$ respectively .

We are going to prove that $\NeighbrAA=\NeighbrBB$ and the map $F$ actually is the identity map.
Indeed, let $\fvar\in\NeighbrAA$ be arbitrary, $p=\pi(\fvar)$ and $\sigma_0=\EqBase(\fvar)$.
Then $\fvar=\sigma_0\RFormBase_p$.
The maps $f_1(z)=\holpar_p(z)$ and $f_2(z)=F\bigl(\holpar_p(z)\bigr)$ are both holomorphic and coincide when $z$ is real.
By analytic continuation $f_1=f_2$.
In particular, we have
$$
	F(\fvar)=F\bigl(\holpar_p(i\sigma_0)\bigr)=\holpar_p(i\sigma_0)=\fvar,
$$
that is the map $F$ is the identity map on $\NeighbrAA$, as claimed.

Now a standard argument implies that there exist a unique real analytic
complex structure $J$ on the set $\domJ$ defined as
the union of all the open neighbourhood $\domJ'$ of $\RManBase$
which admit a complex structure $J'$ which make the pair $(\domJ',J')$
an adapted complex tube on $\RManBase^\RFormBase$.

The pair $(\domJ,J)$ is then the required unique real analytic maximal adapted
complex tube for $\RManBase^\RFormBase$.

This concludes the proof of Theorem \ref{thm::AdaptedRealAnalytic}.

\medskip
{\bf Remark 3.1. }
The construction of the map $\DiffeomB$ is similar to a
one given in \cite{article:DuchampKalka}.

%\end{remark}
%\eject

\end{document}